\documentclass{amsart}

\usepackage{amsfonts,amssymb,mathrsfs}
\usepackage{graphics,color}

\newtheorem{corollary}{Corollary}
\newtheorem{proposition}{Proposition} 
\newtheorem{lemma}{Lemma} 
\newtheorem{theorem}{Theorem} 
\newtheorem{problem}{Problem} 
\theoremstyle{remark}
 
\title[Reachability problems for products of matrices]{Reachability problems for products of matrices in semirings}
\date{February 26, 2003, updated September 29, 2003.}
\author{St\'ephane Gaubert}
\address{INRIA, Domaine de Voluceau, 78153, Le Chesnay C\'edex, France}
\email{Stephane.Gaubert@inria.fr}
\author{Ricardo Katz}
\address{CONICET. Postal address: Dep. of Mathematics, Universidad Nacional de Rosario, Avenida
Pellegrini 250, 2000 Rosario, Argentina}
\email{rkatz@fceia.unr.edu.ar}

\DeclareMathAlphabet{\mathbbold}{U}{bbold}{m}{n}
\def\sreach(#1,#2){\text{\small\sc SReach}(#1,#2)}
\def\srea(#1){\text{\small\sc SReach}(#1)}
\def\mreach(#1,#2){\text{\small\sc MReach}(#1,#2)}
\def\mrea(#1){\text{\small\sc MReach}(#1)}
\def\vreach(#1,#2){\text{\small\sc VReach}(#1,#2)}
\def\vrea(#1){\text{\small\sc VReach}(#1)}
\def\creach(#1,#2){\text{\small\sc CReach}(#1,#2)}
\def\creachp(#1,#2){\text{\small\sc CReach}'(#1,#2)}
\def\sreachp(#1,#2){\text{\small\sc SReach}'(#1,#2)}
\def\mreachp(#1,#2){\text{\small\sc MReach}'(#1,#2)}
\def\vreachp(#1,#2){\text{\small\sc VReach}'(#1,#2)}
\newcommand{\mui}{\mu_{i}}
\newcommand{\reducesto}{\rightarrow}
\newcommand{\reducesfrom}{\leftarrow}
\newcommand{\equivalent}{\leftrightarrow}
\newcommand{\eme}[1]{\bar{#1}}
\newcommand{\emb}[1]{\hat{#1}}
\newcommand{\comp}{\circ}
\newcommand{\zero}{\mathbbold{0}}
\newcommand{\unit}{\mathbbold{1}}
\newcommand{\NEW}[1]{{\em #1}}
\newcommand{\set}[2]{\{#1\mid\,#2\}}
\newcommand{\R}{\mathbb{R}}
\newcommand{\N}{\mathbb{N}}
\newcommand{\Z}{\mathbb{Z}}

\newcommand{\Rm}{\R\cup\{-\infty\}}
\newcommand{\Zp}{\Z\cup\{+\infty\}}
\newcommand{\Np}{\N\cup\{+\infty\}}
\newcommand{\Npm}{\N\cup\{\pm\infty\}}
\newcommand{\Nm}{\N\cup\{-\infty\}}

\newcommand{\nmin}{\N_{\min}}
\newcommand{\nmax}{\N_{\max}}
\newcommand{\nmaxb}{\bar{\N}_{\max}}
\newcommand{\zmin}{\Z_{\min}}
\newcommand{\zmax}{\Z_{\max}}

\newcommand{\nbar}{\bar{\N}}
\newcommand{\sS}{\mathscr{S}}
\newcommand{\sD}{\mathscr{D}}
\newcommand{\sR}{\mathscr{R}}

\newcommand{\sA}{\mathscr{A}}

\newcommand{\sL}{\mathscr{L}}

\newcommand{\sfM}{\mathsf{M}}
\newcommand{\sfF}{\mathsf{F}}

\newcommand{\mrm}[1]{\text{\rm #1}}
\newcommand{\diag}{\mrm{diag}}
\newcommand{\rat}{\mrm{rat}}
\newcommand{\SER}[2]{#1^{\rat}\langle\langle#2\rangle\rangle}
\newcommand{\ser}[1]{\SER{#1}{\Sigma_r}}

\newcommand{\se}{\ser{\sS}}
\newcommand{\npost}{n_{\mrm{P}}}
\newcommand{\nmo}{n_{\mrm{mo}}}
\newcommand{\norm}[1]{\mathsf{n}(#1)}
\newcommand{\vect}{\mrm{vec}}
\begin{document}
\begin{abstract}
We consider the following matrix reachability problem: 
given $r$ square matrices with entries in a semiring, 
is there a product of these matrices which  attains a prescribed matrix? 
We define similarly the vector (resp. scalar)
reachability problem, by requiring 
that the matrix product, acting by right multiplication on a 
prescribed row vector, gives another prescribed row vector
(resp. when multiplied
at left and right by prescribed row and column
vectors, gives a prescribed scalar). 
We show that over any semiring, scalar reachability
reduces to vector reachability which is equivalent to matrix
reachability, and that for any of these problems,
the specialization to any $r\geq 2$ is equivalent
to the specialization to $r=2$.
As an application of this result and of 
a theorem of Krob, we show
that when $r=2$, the vector and matrix reachability
problems are undecidable over the max-plus semiring
$(\Z\cup\{-\infty\},\max,+)$.
We also show that the matrix, vector, and scalar
reachability problems are decidable 
over semirings whose elements
are ``positive'',
like the tropical semiring $(\Np,\min,+)$.
\end{abstract}
\maketitle
\section{Introduction and statement of results}
We consider the following problem:
\begin{problem}[Matrix reachability]\label{pb1}
Given $n\times n$ matrices $A_1,\ldots,A_r$ and $M$
with entries in a semiring $\sS$,
is there a finite sequence $1\leq i_1,\ldots,i_k\leq r$
such that $A_{i_1}\cdots A_{i_k}=M$?
\end{problem}
(Let us recall that a semiring is a set $\sS$ equipped
with an addition and a multiplication,
such that: $\sS$ is a commutative monoid for addition,
$\sS$ is a monoid for multiplication, 
multiplication left and right distributes
over addition, and the zero element for addition
is left and right absorbing for multiplication.)
The matrix reachability problem, which asks whether
$M$ belongs to the semigroup generated by $A_1,\ldots,A_r$,
may be called more classically the {\em semigroup membership} problem.
We chose our terminology to show the interplay
with the two following problems:
\begin{problem}[Vector reachability]\label{pb3} 
Given $n\times n$ matrices $A_1,\ldots,A_r$
and two $1\times n$ matrices $\alpha,\eta$,
all with entries in a semiring $\sS$,
is there a finite sequence $1\leq i_1,\ldots,i_k\leq r$
such that $\alpha A_{i_1}\cdots A_{i_k}=\eta$? 
\end{problem} 
\begin{problem}[Scalar reachability]\label{pb2} 
Given $n\times n$ matrices $A_1,\ldots,A_r$, 
a $1\times n$ matrix $\alpha$, a $n\times 1$ matrix $\beta$, 
all with entries in a semiring $\sS$, and a scalar $\gamma\in\sS$, 
is there a finite sequence $1\leq i_1,\ldots,i_k\leq r$
such that $\alpha A_{i_1}\cdots A_{i_k}\beta=\gamma$? 
\end{problem} 
When $M$ is the zero matrix, the matrix reachability
problem is the well studied {\em mortality problem}. 
Paterson~\cite{paterson70} proved that
when $\sS=(\Z,+,\times)$ is the ring of integers,
the mortality problem is undecidable, even
when $n=3$ and $r=2\npost +2$,
where $\npost$ is the minimal number of pairs of words for which
Post's correspondence problem is undecidable
(Matiyasevich and S\'enizergues~\cite{senizergues} 
proved that $\npost\leq 7$).
Bournez and Branicky~\cite[Prop.~1]{BB02} proved
that the mortality problem remains undecidable
when $n=3$ and $r=\npost+2$,
and Halava and Harju~\cite{HalavaHarju}
proved that the mortality problem
remains undecidable even when $n=3$ and $r=\npost+1$.
See Harju and Karhum\"aki~\cite{karhumaki},
Blondel and Tsitsiklis~\cite{blondel00},
and Halava and Harju~\cite{HalavaHarju} for overviews. 
Blondel and Tsitsiklis~\cite{blondel}, and independently,
Cassaigne and Karhumaki~\cite{cassaigne}, proved that
the mortality problem for $r$ matrices
of dimension $n$ reduces to the mortality problem
for $2$ matrices of dimension $nr$,
which implies that there is an integer 
$\nmo$ such that the mortality problem for two matrices of
dimension $\nmo$ is undecidable
(it follows from~\cite{HalavaHarju} that one can take $\nmo=3(\npost+1)$).

The scalar reachability problem 
previously appeared in the literature in the following form:
\begin{problem}[Corner reachability]\label{pbcorner}
Given $n\times n$ matrices $A_1,\ldots,A_r$
with entries in $\sS$, and a scalar $\gamma\in\sS$,
is there a finite sequence $1\leq i_1,\ldots,i_k\leq r$
such that $(A_{i_1}\cdots A_{i_k})_{1n}=\gamma$? 
\end{problem}
When $\gamma$ is zero,
this becomes the {\em zero corner problem}~\cite{manna,karhumaki,cassaigne},
which is undecidable over $(\Z,+,\times)$ 
when $n=3$ and $r=\npost$~\cite{manna}, and also when 
$r=2$ and $n=3\npost+3$~\cite[Theorem~2 and \S~2.3]{cassaigne}.
An easy observation (\eqref{easy} below) shows
that the scalar and corner reachability problems
are essentially equivalent. 

In this paper, we will show that over any semiring, 
matrix reachability is a problem equivalent to vector reachability,
which is harder than scalar reachability,
and we will also show that for $r\geq 2$,
the $r$-generators version of each of these
problems is equivalent to its $2$-generators variant.
To formalize what ``harder'' and ``equivalent'' means, we have
to define the notion of reduction.
We shall assume that the elements of the semiring $\sS$ 
are represented in some effective way,
and that we have oracles 
taking the representations of two elements $a,b\in \sS$ as input
and returning representations of the sum of $a$ and $b$,
of the product of  $a$ and $b$,
together with the truth value $a=b$, as output.
Then, we say that a problem $P'$
{\em reduces} to a problem $P$, 
and we write $P'\to P$,
if there is an algorithm solving problem $P'$,
using an oracle solving Problem $P$
together with the oracles computing the sum,
the product, and checking the equality 
in $\sS$. The notion we just defined
is a special case of 
{\em Turing reduction}~\cite{hopcroft} with respect to oracles.
We shall also say that $P$ and $P'$ are {\em equivalent},
and write $P\equivalent P'$,
if $P$ reduces to $P'$ and $P'$ reduces
to $P$.

To state more precise results, we need to introduce
restricted versions of the above problems. 
Thus, $\mreach(r,n)$ will denote the specialization of the matrix
reachability problem to $r$ generators of dimension
$n$ and $\mrea(r)$ will denote the specialization
of the matrix reachability problem to $r$ generators
(of arbitrary dimension).
We will use a similar notation
for the vector, scalar, and corner reachability
problems, whose $r,n$ specializations
will be denoted by $\vreach(r,n)$, $\sreach(r,n)$, $\creach(r,n)$,
etc. The following theorem is proved in \S\ref{sec-thmain}.
\begin{theorem}\label{th-main}
In an arbitrary semiring,
\begin{align}
\sreach(r,n)&\mrm{ reduces to } \sreach(2,rn)\label{sreachtosreach}\\
\vreach(r,n)&\mrm{ reduces to } \vreach(2,rn)\label{vreachtovreach}  \\
\mreach(r,n)&\mrm{ reduces to } \mreach(2,rn)\label{mreachtomreach}
\end{align}
\begin{align}
\vreach(r,n)&\mrm{ reduces to }\mreach(r+1,k) \label{vreachtomreach}\\
\nonumber&
\mrm{where }k=n+1\mrm{ if }\eta\neq\zero,\mrm{ and }k=n+3\mrm{ otherwise}
\\
\sreach(r,n)&\mrm{ reduces to }\vreach(r+1,k)\label{sreachtovreach}\\
\nonumber & 
\mrm{where }k=n+1\mrm{ if }\gamma\neq\zero,\mrm{ and }k=n+3\mrm{ otherwise}
\\
\sreach(r,n)&\mrm{ reduces to }\mreach(r+1,k)\label{sreachtomreach}\\
\nonumber & 
\mrm{where }k=n+2\mrm{ if }\gamma\neq\zero,\mrm{ and }k=n+5\mrm{ otherwise.}
\end{align}
Moreover, the value of $\gamma$ 
is preserved in Reduction~\eqref{sreachtosreach}, whereas in Reductions~\eqref{vreachtovreach} and~\eqref{mreachtomreach}, the zero or non zero
character of $\eta$ or $M$ is preserved.
\end{theorem}
Additionnaly, V.~Blondel~\cite{blondel02} observed
that in an arbitrary semiring, 
\begin{align}
\mreach(r,n)\mrm{ reduces to } \vreach(r,rn)
\enspace .\label{mreachtovreach}
\end{align}
For completeness, we reproduce the (simple) proof in~\S\ref{sec-blondel}. 
As an immediate corollary
of Theorem~\ref{th-main} and Reduction~\eqref{mreachtovreach},
we get:
\begin{corollary}\label{th-main-p2}
In an arbitrary semiring,
for all $r,r',r'',r'''\geq 2$,
the scalar reachability problem for $r$ matrices 
is equivalent to the scalar reachability problem
for $r'$ matrices, 
which reduces to the vector reachability problem for $r''$ matrices, which
is equivalent to the matrix reachability problem
for $r'''$ matrices.
\end{corollary}

All the reductions in the proofs of the present paper
take a polynomial time (but the problems
should not be expected to be polynomial, except
in very special cases). 

It would be surprising if the reduction
$\srea(r)\reducesto \mrea(r+1)$ 
stated in~\eqref{sreachtomreach}
could be improved to give $\srea(r)\reducesto\mrea(r)$.
Indeed, when $\sS=\Z$, 
$\srea(1)$ is equivalent to the 
Pisot problem, a well known unsolved problem
consisting in deciding the existence of a zero in
an integer linear recurrent sequence,
whereas the matrix reachability
problem $\mrea(1)$ becomes:
\begin{align*}
\mrm{given }A,M\in\Z^{n \times n},
\mrm{ is there some }k\geq 1\mrm{ such that }A^k=M \enspace,
\end{align*}
a much simpler problem, which is 
even solvable in polynomial time, see~\cite{kannan},
and also~\cite{cai,babai}.
The vector reachability problem for one matrix, $\vrea(1)$,
which was called the {\em orbit problem} in~\cite{kannan},
is also solvable in polynomial time~\cite{kannan},
so that the existence of a reduction $\srea(r)\reducesto\vrea(r)$
seems unlikely. 

In the statement of Theorem~\ref{th-main},
we needed to distinguish the cases where $M$, $\eta$, or $\gamma$, are zero.
Indeed, the reductions depend critically on the zero or non-zero character
of the instance. For instance, the proof of~\eqref{mreachtomreach}
when $M=\zero$ follows merely from the argument
of Blondel and Tsitsiklis~\cite{blondel} and of
Cassaigne and Karhumaki~\cite{cassaigne}, 
whereas the $M\neq\zero$ case is proved 
using a very different method (compare \S\ref{subsec-btck}
with \S\ref{subsec-ournew}). 

We next derive some consequences of Theorem~\ref{th-main}.
Let us consider the case when $\sS$ is the max-plus semiring 
$\zmax=(\Z\cup\{-\infty\},\max,+)$. In $\zmax$, 
the matrix product is given by 
\[ 
(AB)_{ij}=\max_{k} (A_{ik}+B_{kj})\enspace. 
\] 
In $\zmax$, the scalar reachability problem was solved negatively
by Krob:
\begin{theorem}[See~\cite{krob93b}]\label{th-krob} 
For $r=2$, the scalar reachability problem 
over the max-plus semiring $\zmax$ is undecidable.  
\end{theorem} 
In fact, Krob did not make explicit Theorem~\ref{th-krob}, 
but we shall see in \S~\ref{sec-krob} that Theorem~\ref{th-krob} 
is contained in his proof.
Note also that Krob stated
the results in the {\em equatorial semiring} $\zmin=(\Zp,\min,+)$,
which is effectively isomorphic to $\zmax$ (by a change of
sign), so that decidability issues in $\zmax$ and $\zmin$
are equivalent. 
We get as a corollary of Theorem~\ref{th-krob}
and of the reductions~\eqref{sreachtovreach}, \eqref{sreachtomreach}, 
\eqref{vreachtovreach}, and \eqref{mreachtomreach} in Theorem~\ref{th-main}:
\begin{theorem}\label{th-1} 
For $r=2$, the matrix and vector reachability problems 
over the max-plus 
semiring $\zmax$ are undecidable.
\end{theorem} 
The $r\geq 2$ bound
is optimal, since when $r=1$, the matrix reachability problem
in $\zmax$ is known to be decidable
(see~\S\ref{sec-r=1} below).
Moreover, a simple argument shows that for any $r$,
the {\em mortality} problem 
in the max-plus semiring is decidable
(use the third remark after Theorem~2 in~\cite{blondel}). 

The proof of Theorem~\ref{th-1} and Krob's proof of Theorem~\ref{th-krob}
show that the restrictions of the scalar, vector, and matrix
reachability problems to matrices of some fixed, sufficiently large,
dimension $n$, remain undecidable.
Indeed, Matiyasevich's theorem (see in particular the corollary
in the introduction of~\cite{dmr76}, \cite{mbook}, and the references therein)
shows that the Hilbert's tenth problem remains undecidable
for a subclass of instances consisting of a family
of polynomials of bounded degree, with a fixed number
of variables, and one can check that Krob's proof,
when applied to this family, yields linear representations
of bounded dimension $n$. 

A natural question
would be to find an alternative proof which would allow
a more precise control of the dimension. 
The reader should note, here, 
that the Post-correspondence based
technique of Paterson~\cite{paterson70},
which relies on the embedding
of a free monoid with at least two letters into 
matrices over $\Z$, has no natural extension
to $\zmax$. In fact, the possibity of such an extension
was considered when the equality problem for max-plus rational
series was still open, and it was remarked
independently by Krob and by Simon~\cite{simon88},
that $\zmax^{n\times n}$ contains no free submonoid. 
(To see this, define, for all $A\in \zmax^{n\times n}$, 
$\norm A = \sup\set{|A_{ij}|}{1\leq i,j\leq n,\;A_{ij}\neq -\infty}$,
observe that $\norm{AB}\leq \norm A + \norm B$,
and deduce that any finitely generated matrix submonoid of $\zmax^{n\times n}$
has a growth function $O(k^{n^2})$.)

As easy corollaries of Theorem~\ref{th-1},
we will get in \S\ref{sec-proj} 
undecidability results for {\em projective}
variants of the reachability problem.
Recall that the proportionality
relation $\sim$ on $\zmax^n$ and $\zmax^{n\times n}$
is defined by $u\sim v$ if $u=\lambda v$, for
some $\lambda \in \R$ (that is, $u_i=\lambda+ v_i$
when $u,v\in \zmax^n$, or $u_{ij}=\lambda+ v_{ij}$
when $u,v\in \zmax^{n\times n}$). 
\begin{corollary}[Projective matrix reachability over $\zmax$ is undecidable] \label{cor-proj}
The following problem is undecidable:
given $A_1,A_2,M\in \zmax^{n\times n}$,
is there a finite sequence $1\leq i_1,\ldots,i_k\leq 2$
such that $A_{i_1}\cdots A_{i_k}\sim M$?
\end{corollary}
\begin{corollary}[Projective vector reachability over $\zmax$ is undecidable] \label{cor-proj2}
The following problem is undecidable:
given  $A_1,A_2\in \zmax^{n\times n}$,
and $\alpha,\eta\in\zmax^{1\times n}$, 
is there a finite sequence $1\leq i_1,\ldots,i_k\leq 2$
such that $\alpha A_{i_1}\cdots A_{i_k}\sim \eta$?
\end{corollary}
Projective reachability problems
arise in relation with the problem 
of determining whether a max-plus rational
series is subsequential (i.e. has a deterministic linear representation).
A general result which was first understood by 
Choffrut, see~\cite[Ch.~3]{choffrut} and~\cite[Th.~1]{choffrut01} 
for a recent overview, see also~\cite[Th.~4]{gaubert95c}
and~\cite[Th.~10]{mohri}, yields a partial-decision algorithm to 
determine whether a max-plus 
rational series is subsequential
(by ``partial decision'', we mean that the algorithm need not terminate, even
when the series is subsequential).
In the case of $\zmax$, this algorithm 
consists in computing the set
$\set{\alpha A_{i_1}\ldots A_{i_k}}{k\geq 1, \; 1\leq i_1,\ldots,i_k\leq r}$
modulo the equivalence relation $\sim$.
Thus, Corollary~\ref{cor-proj2} shows that it is 
undecidable whether this
algorithm will produce the equivalence class of $\eta$.
Corollary~\ref{cor-proj2} also has an interesting
discrete event systems interpretation.
In this context~\cite{gaubert95c,gaubert96a,gaumair95,CAV,gaumair97a},
the vector $\alpha':=\alpha A_{i_1}\ldots A_{i_k}$ gives
the completion time of different events, after
the execution of a schedule represented by a
sequence $i_1,\ldots,i_k$, and the equivalence
class of $\alpha'$ modulo $\sim$
represents inter-event delays
($\alpha'_i-\alpha'_j$ represents typically
the time a part stays in a storage resource).

Even, in the classical case $\sS=\Z$,
Theorem~\ref{th-main} suggests some results.
For instance, the result of~\cite{manna,karhumaki,cassaigne}
showing that the zero corner problem
is undecidable over $\Z$ implies that
the vector reachability problem
is undecidable, and this does not
seem to have been stated previously. In fact, we can
prove a more precise result by a small modification of the proof
of Paterson~\cite{paterson70}
(see \S\ref{sec-proof-cor-vmort}):
\begin{proposition}\label{cor-vmort}
The vector mortality problem over $\Z$,
for $\npost+1$ matrices of dimension $3$,
is undecidable. 
\end{proposition}
Theorem~\ref{th-main} also improves
by one unit the dimension obtained
in Theorem~2 of  \cite{cassaigne}
(see \S\ref{sec-cor-cassaigne}):
\begin{corollary}\label{cor-cassaigne}
The zero corner problem 
for $2$ matrices of dimension $3\npost+2$
over $\Z$ is undecidable. 
\end{corollary}
It is natural to ask whether the reachability problems become
decidable in other semirings.
Many variants of $\zmax$ can be found in the literature. In particular, the following semirings are listed in~\cite{pin95}: 
\begin{center} 
\begin{tabular}{ccc} 
$\nmin = (\Np,\min,+)$ & Tropical semiring~\cite{simon90}\\ 
$\nmax = (\Nm,\max,+)$ & Boreal semiring~\cite{krob93b}\\ 
$\nmaxb = (\Npm,\max,+)$ &  Mascle's semiring~\cite{mascle}\\ 
$\sL=(\N\cup\{\omega,+\infty\},\min,+)$ & Leung's semiring~\cite{leung91} 
\end{tabular} 
\end{center} 
In $\nmaxb$, $(+\infty)+(-\infty)=(-\infty)+(+\infty)=-\infty$. 
Leung's semiring $\sL$ is the one point compactification of the semiring $\nmin$ equipped with its 
discrete topology: the minimum is defined with respect to the order $0<1<2<\cdots <\omega<+\infty$, 
and the addition of $\nmin$ is completed by $\omega+a=a+\omega=\max(a,\omega)$. 
The semiring $\nmax$ is a subsemiring of $\zmax$, and the map $x\mapsto -x $ is an isomorphism from $\nmin$ 
to a subsemiring of $\zmax$. 

To show that the reachability problems are decidable over these semirings 
we will need the following definitions. We say that a semiring
$\sS$ {\em is separated by morphisms of finite image} if for all 
$\gamma \in \sS$, 
there is a finite semiring $\sS_\gamma $ and a semiring morphism $\pi_\gamma $ from $\sS$ to $\sS_\gamma $ 
such that $\pi_{\gamma }^{-1}(\pi_{\gamma }(\gamma ))=\{\gamma \}$. 
We shall say that $\sS$ is {\em effectively} separated by morphisms
of finite image if the maps $\gamma\mapsto \sS_\gamma$ and 
$\gamma\mapsto \pi_{\gamma}$
are effective, in the sense that for any $\gamma \in \sS$,
we can compute the (finite) set of elements of
the semiring $\sS_\gamma$,
together with the addition and multiplication tables of $\sS_\gamma$,
and that we can compute $\pi_\gamma (y)$ for any 
$y\in \sS$. We prove in~\S\ref{sec-decidable}:
\begin{theorem}\label{th-4} 
The matrix, vector, and scalar reachability
problems are decidable over a semiring
that is effectively separated by morphisms
of finite image.
\end{theorem} 
We will show in fact a slightly more precise result
(Theorem~\ref{th-5} in~\S\ref{sec-decidable}).

Our method applies not only to max-plus type semirings, but 
also to the semiring of natural numbers,
$\N = (\N,+,\times)$, and to its
completion, $\nbar = (\Np,+,\times)$
(in $\nbar$, we adopt the convention
$0\times (+\infty) =(+\infty)\times 0=0$).
A simple argument, which is given in \S~\ref{sec-proof-prop-1}, shows
that:
\begin{proposition}\label{prop-1}
The semirings $\nmin $, $\nmax $, $\nmaxb $, $\sL$, $\N$, and $\nbar$,
all are effectively separated by morphisms of finite image.
\end{proposition}
As a corollary of Theorem~\ref{th-4} and Proposition~\ref{prop-1} we get:
\begin{corollary}\label{cor-dec}
The matrix, vector, and scalar reachability problems are decidable over the semirings $\nmin$, $\nmax,\nmaxb$, $\sL$, $\N$, and $\nbar$. 
\end{corollary}
When the semiring is $\nmin$ or $\nmax$, the
decidability of the scalar reachability
problem was stated by Krob in~\cite[Proposition~2.2]{krobfatou}.

We already observed from Theorem~\ref{th-main}
that matrix reachability, or equivalently,
vector reachability, are harder problems than
scalar reachability.
However, for all the examples of semirings that we considered,
either all problems were undecidable, or they were all decidable.
This raises the question of the existence
of a semiring with undecidable
matrix reachability problem but decidable scalar reachability
problem. 

Let us finally mention some additional 
motivation and point out other references. 
Automata with multiplicities and semigroups of matrices over the 
tropical semiring have been much studied
in connection with decision problems in language theory,
see~\cite{simon78,simon94}, \cite{hashiguchi,hashiguchi90}, \cite{mascle}, \cite{leung91}, \cite{krob93b}, and~\cite{pin95} for a survey.
Automata with multiplicities over the max-plus semiring and
max-plus linear semigroups appear in the modelling 
of discrete event dynamic systems,
see~\cite{gaubert95c,gaumair95,CAV,gaumair97a,gaubert99f}.
General references about max-plus algebra 
are~\cite{bcoq,cuning,maslovkolokoltsov95,gondran02}. 
Some of the present results have been announced in~\cite{gk03}.

\medskip\noindent{\em Acknowledgments.}\/
The authors thank Vincent Blondel and Daniel Krob 
for helpful comments on a preliminary version of
this manuscript.

\section{Preparation}\label{sec-proof}
In this section, we collect preliminary results.
\subsection{Reformulation in terms of linear representations}
The proof of the results uses 
rational series and automata notions: we next recall basic 
definitions. See~\cite{berstelreut} or~\cite{lallement}
for more background.

Let $\Sigma_r=\{a_1,\ldots,a_r\}$ denote an alphabet with $r$ letters,
 and let $\Sigma_r^*$ denote the \NEW{free monoid} 
on $\Sigma_r$, that is the set of finite (possibly empty) words
with letters in $\Sigma_r$. A subset of $\Sigma_r^*$ is called 
a \NEW{language}. We will also consider
the {\em free semigroup} $\Sigma_r^+$, which
is the subsemigroup of $\Sigma_r^*$ composed
of nonempty words.
We say that a map $s:\Sigma_r^*\to \sS$ is {\em recognizable} or {\em rational} if there exists an integer $n$,  $\alpha\in \sS^{1\times n}$, $\beta\in \sS^{n\times 1}$, and a morphism $\mu: \Sigma_r^*\to \sS^{n\times n}$ 
such that $s(w)=\alpha\mu(w)\beta$ for all $w\in \Sigma_r^*$. We say that $(\alpha,\mu,\beta)$ is 
a \NEW{linear representation} of $s$, and that $n$ is the \NEW{dimension} of the representation. We denote by $\ser\sS$ 
the set of rational maps, which are also called \NEW{rational series}.
Problems~\ref{pb1}--\ref{pbcorner}, 
can be rewritten as: 
\begin{align}
\label{def-mreach}
&\mreach(r,n):\\
&\qquad \mu\mrm{ morphism } \Sigma_r^*\to \sS^{n\times n},M\in \sS^{n\times n};\; \exists w\in \Sigma_r^+, \mu(w)=M \enspace ?
\nonumber
\\
\label{def-vreach}&\vreach(r,n):\\
&\qquad \mu\mrm{ morphism } \Sigma_r^*\to \sS^{n\times n},\alpha,\eta\in \sS^{1\times n};\; \exists w\in \Sigma_r^+, \alpha\mu(w)=\eta \enspace ?\nonumber
\\
\label{def-sreach} &\sreach(r,n):\\
\nonumber &\qquad s\in \ser{\sS}
\mrm{ with a linear representation of dimension }n, 
\gamma\in\sS;\\  &\qquad \exists w\in \Sigma_r^+, s(w)=\gamma \enspace ?
\nonumber \\
\label{def-creach}&\creach(r,n):\\
&\qquad \mu\mrm{ morphism } \Sigma_r^*\to \sS^{n\times n},\gamma\in \sS;\; 
\exists w\in \Sigma_r^+, \mu_{1,n}(w)=\gamma \enspace ?\nonumber
\end{align} 
\subsection{Variants allowing the empty word} 
We required that the word $w$ belongs to $\Sigma_r^+$ and not
to $\Sigma_r^*$ in the
formulations~\eqref{def-mreach}--\eqref{def-creach},
because in the statements of 
Problems~\ref{pb1}--\ref{pbcorner},
we considered sequences $i_1,\ldots,i_k$ of length at least $1$.
However, some simple observations will show that
putting $w\in\Sigma_r^*$ or $w\in\Sigma_r^+$ 
in~\eqref{def-mreach}--\eqref{def-creach} is 
essentially irrelevant.

Let us denote by $\mreachp(r,n)$,
$\vreachp(r,n)$, $\sreachp(r,n)$, and $\creachp(r,n)$ the variants
of the above problems with $\Sigma_r^*$ instead
of $\Sigma_r^+$ in~\eqref{def-mreach}--\eqref{def-creach}.
We will denote by $\zero$ and $\unit$ the zero and unit
elements of $\sS$, respectively,
by $\zero_{pq}\in \sS^{p\times q}$
or simply $\zero$ the $p\times q$ zero matrix, and by $I_n\in \sS^{n\times n}$
or simply $I$ the $n\times n$ identity matrix. 
Since the cases where $M=I$ or $M=\zero$, or $\eta=\zero$,
or $\gamma=\zero$ will sometimes require a special
treatment, we will incorporate restrictions about $M$,
$\eta$, or $\gamma$ in the notation,
writing for instance $\mreach(r,n,M\neq I)$ 
for the restriction of the matrix reachability 
problem to $r$ generators of dimension $n$
and a matrix $M$ different from the identity. 
\begin{lemma}
\label{lem-triv}
The following reductions hold:
\begin{align}
\label{equi-vreach}
\vreach(r,n)&\equivalent \vreachp(r,n)
\end{align}
\begin{align}
\label{equi-sreach}
\sreach(r,n)&\equivalent \sreachp(r,n)
\end{align}
\begin{align}
\label{equi-mreach1}
\mreach(r,n,M\neq I)&\equivalent \mreachp(r,n,M\neq I)
\end{align}
\begin{align}
\mreach(r,n,M=I)& \reducesto \mreachp(r,n+1,M\neq I) 
\label{equi-mreach2}
\end{align}
\begin{align}
\creach(r,n\geq 2,\gamma\neq \zero)&\equivalent 
\creachp(r,n\geq 2,\gamma\neq \zero)
\label{equi-creach}
\enspace .
\end{align}
\end{lemma}
\begin{proof}
Consider an instance of $\vreach(r,n)$,
which consists of $\alpha,\eta\in \sS^{1\times n}$
and a morphism $\mu:\Sigma_r^*\to \sS^{n\times n}$.
Since
\begin{align}
(\exists w\in \Sigma_r^+ ,\; \alpha\mu(w)=\eta)
\iff (\exists a\in \Sigma_r,\; \exists z\in \Sigma_r^*,\; \alpha\mu(a) \mu(z)=\eta)
\enspace, 
\end{align}
$\vreach(r,n)$ reduces to $\vreachp(r,n)$.
Conversely, 
\begin{align}
(\exists w\in \Sigma_r^* ,\; \alpha\mu(w)=\eta)
\iff (\alpha=\eta\mrm{ or } \exists w\in \Sigma_r^+,\; \alpha \mu(w)=\eta)
\enspace, 
\end{align}
shows that $\vreachp(r,n)$ reduces to $\vreach(r,n)$,
which shows~\eqref{equi-vreach}.

Similarly, consider an instance
of $\sreach(r,n)$, consisting
of $s\in \ser{\sS}$ with a linear representation
of dimension $n$, $(\alpha,\mu,\beta)$,
and $\gamma\in \sS$. Let $\_$ denote the empty word
of $\Sigma_r^*$, and let $a^{-1}s$
denote the series defined by $a^{-1}s(w):=
s(aw)$. Observe that $a^{-1}s$ is recognized
by the linear representation $(\alpha\mu(a),\mu,\beta)$,
which is is still of dimension $n$.
Since
\[
(\exists w\in \Sigma_r^+, \;s(w)=\gamma)\iff
(\exists a\in \Sigma_r,\; \exists w\in \Sigma_r^*,\;
a^{-1}s(w)= \gamma)\enspace,
\]
$\sreach(r,n)$ reduces to $\sreachp(r,n)$. 
Conversely,
\begin{align}
(\exists w\in \Sigma_r^*,\;s(w)=\gamma)
\iff (s(\_)=\gamma \mrm{ or } \exists w\in \Sigma_r^+,
\;s(w)=\gamma) \enspace,
\end{align}
shows that $\sreachp(r,n)$ reduces to $\sreach(r,n)$,
which shows~\eqref{equi-sreach}. 

The problems $\mreachp(r,n,M\neq I)$ and 
$\mreach(r,n,M\neq I)$
are trivially equivalent, because $\mu$ sends
the empty word to the identity matrix.
This shows~\eqref{equi-mreach1}.

Before showing~\eqref{equi-mreach2},
we introduce a notation that we
shall use repeatedly in the sequel. 
If $U_1,\ldots,U_k$ are square matrices
with entries in a semiring $\sS$,
we denote by $\diag(U_1,\ldots,U_k)$ the
block diagonal matrix whose diagonal blocks
are $U_1,\ldots,U_k$. If $\mu_1,\ldots,\mu_k$
are morphisms from a free monoid to matrix
monoids, we denote by $\diag(\mu_1,\ldots,\mu_k)$
the morphism which sends a word $w$ to $\diag(\mu_1(w),
\ldots,\mu_k(w))$. 
For all $1\leq p$, we denote by 
$\zero_{pp}$ the zero morphism from
a free monoid to $\sS^{p\times p}$.
Let us define now the morphism $\mu':\Sigma_r^*\to\sS^{(n+1)\times (n+1)}$, 
$\mu'=\diag(\mu,\zero_{11})$.
Since $\mu'(\_)= I$, 
\[
(\exists w\in \Sigma_r^+,\; \mu(w)=I)
\iff (\exists w\in \Sigma_r^*,\; \mu'(w)=\diag(I,\zero_{11})) \enspace,
\]
which shows~\eqref{equi-mreach2}.

Finally, the problems $\creachp(r,n\geq 2,\gamma\neq \zero)$ and 
$\creach(r,n\geq 2,\gamma\neq \zero)$
are trivially equivalent, because $\mu(\_)_{1n}=\zero\neq \gamma$,
as soon as $n\geq 2$ and $\gamma\neq \zero$. 
This shows~\eqref{equi-creach}. 
\end{proof}
We did not consider the problems
$\mreachp(r,n,M= I)$ and
$\creachp(r,n\geq 2,\gamma=\zero)$
in Lemma~\ref{lem-triv},
since the answer to these problems
is trivially ``yes''.
\subsection{Equivalence of corner and scalar reachability}
The following elementary reductions show that
up to an increase of the dimension 
of matrices, the corner and scalar reachability
problems are equivalent:
\begin{align}
\label{easy}
\sreach(r,n) \reducesto \creach(r,n+2) \enspace ,\\
\label{easy2}
\creach(r,n) \reducesto \sreach(r,n) \enspace .
\end{align}
Indeed, consider an instance of $\sreach(r,n)$,
which consists of a series $s\in \se$ 
with a linear representation $(\alpha,\mu,\beta)$ of dimension $n$,
and a scalar $\gamma \in \sS$. 
We build the morphism $\mu': \Sigma_r^*\to\sS^{(n+2)\times(n+2)}$ such that 
\begin{align}
\mu'(a_i) = 
\begin{pmatrix} 
\zero_{11} & \alpha\mu(a_i) & \alpha\mu(a_i)\beta\\  
\zero_{n1} & \mu(a_i) & \mu(a_i)\beta \\  
\zero_{11} & \zero_{1n} & \zero_{11} 
\end{pmatrix} 
\qquad\forall 1\leq i\leq r \enspace .
\label{eq-def-mup}
\end{align}
An immediate induction on the length of $w$ shows that 
\begin{align}
\label{eq-p-mup} \mu'(w)= 
\begin{pmatrix} 
\zero_{11} & \alpha\mu(w) & s(w) \\  
\zero_{n1} & \mu(w) & \mu(w)\beta \\  
\zero_{11} & \zero_{1n} & \zero_{11} 
\end{pmatrix} 
\enspace ,\quad \forall w\in\Sigma_{r}^+ \enspace . \end{align} 
Thus,
\begin{align}
\label{eq-scaltocorn}
\forall w\in \Sigma_r^+,\; 
\mu'_{1,n+2}(w)=s(w)
\enspace ,
\end{align}
which shows~\eqref{easy}.

Reduction~\eqref{easy2} holds
because $\creach(r,n)$ is merely
a special case of $\sreach(r,n)$. 
Indeed, consider an instance of the corner
reachability problem,
consisting of $\mu$ as above, and $\gamma\in \sS$,
and let $\alpha=(\unit,\zero,\ldots,\zero)\in \sS^{1\times n}$,
and $\beta=(\zero,\ldots,\zero,\unit)^T\in \sS^{n\times 1}$.
Then, for all $w\in \Sigma_r^*$,
$\mu_{1n}(w)=\alpha \mu(w)\beta$, 
which shows~\eqref{easy2}. 
\subsection{Matrix representation of trim unambiguous automata}
\label{sec-cons-mat}
We shall use several times the following essentially classical
constructions. To any automaton $\sA$
over $\Sigma_r$ with set of states $\{1,\ldots,p\}$,
and set of initial (resp. final) states $I$ (resp. $F$),
we associate the morphism
$\nu_{\sA}: \Sigma_r^*\to\sS^{p\times p}$,
\begin{align}
\label{e-def-auto}
\forall x\in \Sigma_r,\quad \nu_{\sA}(x)_{ij} = 
\begin{cases} \unit & \mrm{if there is an arrow from } i\mrm{ to }j\mrm{ labeled }x\mrm{ in }\sA,\\ 
\zero & \mrm{otherwise,} 
\end{cases} \end{align}
the vectors
\[
\alpha_{\sA}\in \sS^{1\times p}, \; (\alpha_{\sA})_k=
\begin{cases}
\unit & \mrm{if } k\in I\\
\zero & \mrm{otherwise,}
\end{cases}
\quad 
\beta_{\sA}\in \sS^{p\times 1}, \; (\beta_{\sA})_k=
\begin{cases}
\unit & \mrm{if } k\in F\\
\zero & \mrm{otherwise,}
\end{cases}
\]
together with 
\begin{align*}
\sfM_{\sA}&=\set{\nu_{\sA}(v)}{v\in \Sigma_r^* \mrm{ and } [\nu_{\sA}(v)]_{\iota \phi}=\unit \mrm{ for some }\iota\in I,\;\phi\in F} \enspace ,\\
\sfF_{\sA}&=\set{\alpha_{\sA}\nu_{\sA}(v)}{v\in \Sigma_r^* \mrm{ and }
 [\alpha_{\sA}\nu_{\sA}(v)]_{\phi}=\unit \mrm{ for some }\phi\in F} \enspace .
\end{align*}
Recall that $\sA$ is {\em unambiguous} if for all
$w\in \Sigma_r^*$, there is at most one path with label
$w$ from an input state to an output state,
and that $\sA$ is {\em trim} if for all state $k$,
there is a path from some input state to $k$,
and a path from $k$ to some output state.
\begin{lemma}\label{lem-cons}
If $\sA$ is trim and unambiguous, then $\sfM_{\sA}$ and $\sfF_{\sA}$
can be effectively
computed, and the language $L$ recognized by $\sA$
is 
$\set{w\in \Sigma_r^*}{\nu_{\sA}(w)\in \sfM_\sA}=
\set{w\in \Sigma_r^*}{\alpha_{\sA}\nu_{\sA}(w)\in \sfF_{\sA}}$.
\end{lemma}
\begin{proof}
If $\sA$ is trim and unambiguous, for all $1\leq i,j\leq p$
and $w\in \Sigma_r^*$, there is at most one path from $i$
to $j$ with label $w$. 
Then, it follows from the well known graph
interpretation of the matrix product (see e.g.~\cite[\S~4.7]{stanley}),
that all the matrices $\nu_{\sA}(v)$ have $\zero,\unit$ entries
(which implies that $\sfM_{\sA}$ is finite and can be effectively computed),
and that $L=\nu_{\sA}^{-1}(\sfM_{\sA})$. The
analogous property for $\sfF_{\sA}$ 
is proved in a similar way.
\end{proof}

\subsection{Derivation of Theorem~\protect\ref{th-krob} from Krob's proof}
\label{sec-krob} 
Krob considered the following problems for 
series $s,t\in \ser{\sS}$ and $\sS=\nmin,\nmax,\zmin$:
\begin{eqnarray*} 
(\mrm{Equality})&& s,t\in \ser{\sS};\;s=t\enspace?\\ 
(\mrm{Inequality})&& s,t\in \ser{\sS};\;s\leq t\enspace?\\ 
(\mrm{Local Inequality})&& s,t\in \ser{\sS};\; 
\exists w\in \Sigma_r^*,\; s(w) \leq t(w)\enspace?\\ 
(\mrm{Local Equality})&& s,t\in \ser{\sS};\; 
\exists w\in \Sigma_r^*,\; s(w) = t(w)\enspace?  
\end{eqnarray*} 
Corollary~4.3 of~\cite{krob93b} shows that all these problems are undecidable when $\sS=\nmin$ 
or $\sS=\nmax$, provided that the number of letters $r$ is at least $2$. The undecidability of the scalar reachabiliy problem does not follow from this 
statement, but it does follow from the proof of~\cite{krob93b}. 
Indeed, in \S~3 of~\cite{krob93b}, 
Krob associates effectively to any instance $(I)$ of Hilbert's tenth problem a rational series denoted by $HD$, with coefficients 
in $\zmin$, over an alphabet $A$, with the property that $HD(w)\leq 0$ for all $w\in A^*$ 
and that there is a word $z\in A^*$ such that $HD(z)=0$ if, and only if, instance $(I)$ has a solution. Since Hilbert's tenth 
problem is undecidable, this implies that the scalar reachability problem over the semiring $\zmin$ is undecidable, when $\gamma=0$. 
Moreover, the coding argument given at the beginning of the proof of Theorem~3.1 of~\cite{krob93b} associates effectively to $HD$ 
a rational series $\sigma(HD)$ with coefficients in $\zmin$ 
over a two letters alphabet, and this series takes the same finite values as $HD$. This shows Theorem~\ref{th-krob}. 
\section{Embedding matrix semigroups with $r$ generators in matrix semigroups
with $2$ generators}
The proof of Theorem~\ref{th-main} relies on two different embeddings
of semigroups of $n\times n$ matrices with $r$-generators 
in semigroups of $nr\times nr$ matrices with  $2$-generators.
\subsection{First embedding}\label{sec-emb1}
Let $b,c$ denote two letters.
To any morphism $\mu:\Sigma_r^*\to \sS^{n\times n}$,
we associate the morphism $\eme\mu:\{b,c\}^*\to \sS^{nr\times nr}$,
defined by:
\begin{align}
\label{eq-def-kappa}
\eme\mu(b)=
\begin{pmatrix} \mu(a_1)& \zero_{n,(r-1)n}\\
\vdots & \vdots\\ 
\mu(a_{r}) & \zero_{n,(r-1)n} 
\end{pmatrix}, 
\quad\mrm{and} \quad \eme\mu(c)=
\begin{pmatrix} 
\zero_{(r-1)n,n}& I_{(r-1)n}\\ 
\zero_{n,n} & \zero_{n,(r-1)n} 
\end{pmatrix} 
\end{align} 
(recall that $I_n$ denotes the $n\times n$ identity matrix). 

Following an usual device, we shall associate
to any word of $\Sigma_r^*$ a word of $\{b,c\}^*$
by way of the coding function $\delta: \Sigma_r^*\to\{b,c\}^*$, 
\begin{align}
\label{e-def-code}
\delta(a_{i_1}\ldots a_{i_k}) = c^{i_1-1} b\ldots c^{i_k-1}b 
\enspace,
\end{align}
for all $1\leq i_1,\ldots,i_k\leq r$. 
The function $\delta$ is a bijection from $\Sigma_r^*$ to the language 
$\delta(\Sigma_r^*)=\{b,cb,\ldots,c^{r-1}b\}^*$. 
The following result can be proved by an
immediate induction on $k$.
\begin{proposition}
For all $a_{i_1},\ldots,a_{i_k}\in \Sigma_r$, 
\begin{align}
\label{theformula}
\eme\mu\comp\delta(a_{i_1}\ldots a_{i_k}) =  
\begin{pmatrix} 
\mu(a_{i_1} a_{i_2}\ldots a_{i_k})& \zero_{n,(r-1)n}\\ 
\mu(a_{i_1+1} a_{i_2}\ldots a_{i_k})& \zero_{n,(r-1)n}\\ 
\vdots& \vdots\\ 
\mu(a_r a_{i_2}\ldots a_{i_k})&\zero_{n,(r-1)n} \\ 
\zero_{(i_1-1)n,n}& \zero_{(i_1-1)n,(r-1)n}  
\end{pmatrix} 
\enspace . \end{align}
\end{proposition}
We shall use in particular the
specialization of~\eqref{theformula} to $i_1=r$:
\begin{align} 
\forall z\in \Sigma_r^*,\;
\eme\mu\comp\delta(a_{r}z) =  
\begin{pmatrix} 
\mu(a_r z)& \zero_{n,(r-1)n}\\ 
\zero_{(r-1)n,n}& \zero_{(r-1)n,(r-1)n} 
\end{pmatrix} 
\enspace. \label{e-news} 
\end{align} 
\subsection{Second embedding}\label{sec-emb2}
This embedding is borrowed
from the proof of~\cite[Th.~1]{blondel}
and~\cite[Th.~1]{cassaigne}.
To any morphism $\mu:\Sigma_r^*\to\sS^{n\times n}$,
we associate the morphism $\emb\mu:\{b,c\}^*\to\sS^{rn\times rn}$:
\begin{align}
\label{cons-blondel}
\emb\mu(b)= \begin{pmatrix} \zero_{(r-1)n,n} & I_{(r-1)n}\\
I_{nn}& \zero_{n,(r-1)n}
\end{pmatrix}, 
\;\;
\emb\mu(c) = \diag(\mu(a_1),\ldots,\mu(a_r))
\enspace .
\end{align}
To simplify notations,
we will use a convention of cyclic indexing
of the letters of $\Sigma_r$, so that $a_{r+1}=a_1$,
$a_{r+2}=a_2$, etc. 
We shall use the trivial fact that any
word $v\in \{b,c\}^*$ can be written
(uniquely) as:
\begin{align}
\label{codefactor}
v=c^{i_1}b\ldots c^{i_k}bc^{i_{k+1}}
\end{align}
where $k\geq 0$ and $0\leq i_1,\ldots,i_{k+1}$,
with the convention that $v=c^{i_{k+1}}=c^{i_1}$
when $k=0$. 
\begin{lemma}\label{lemma-semidirect}
If $v\in \{b,c\}^*$ is written as in~\eqref{codefactor}, 
then, 
\begin{align}
\emb\mu(v)=\diag[ \mu(a_1^{i_1} a_2^{i_2}\ldots a_{k+1}^{i_{k+1}}),
\ldots,
\mu(a_r^{i_1} a_{r+1}^{i_2}\ldots a_{r+k}^{i_{k+1}}) ]\emb\mu(b^k)
\enspace .\label{e-jj}
\end{align}
\end{lemma}
For instance, when $r=3$, \eqref{e-jj} states that:
\[
\emb\mu(c^2bc^7bc^{9}bcb^2c^{11})
=\begin{pmatrix}
\mu(a_1^2a_2^7a_3^9a_1a_3^{11})& \zero & \zero\\
\zero &\mu(a_2^2a_3^7a_1^9a_2a_1^{11}) & \zero\\
\zero & \zero &\mu(a_3^2a_1^7a_2^9a_3a_2^{11}) 
\end{pmatrix}
\emb\mu(b^5)
\enspace .
\]
\begin{proof}[Proof of Lemma~\ref{lemma-semidirect}]
Consider the semigroup $\sD\subset \sS^{rn\times rn}$ 
of block diagonal matrices with $r$ diagonal
blocks of dimension $n$, together with the group
$\sR$ generated by the matrix $B:=\emb\mu(b)$
($B$ is invertible since $B^{-1}=B^T$). 
Since $B\sD B^{-1}\subset \sD$,
for all $(D,R),(D',R')\in \sD\times \sR$,
we have
\begin{align}
DRD'R'= DRD'R^{-1} RR', \mrm{ where } DRD'R^{-1}\in \sD, \; RR'\in \sR \enspace.
\label{semidirect} 
\end{align}
(In other words, the semigroup of matrices of the form $DR$,
where $(D,R)\in \sD\times \sR$, 
is a semidirect product of $\sD$ by $\sR$.)
Then, Formula~\eqref{e-jj} is proved by an
immediate induction, thanks to~\eqref{semidirect},
and to the observation that $D\mapsto B^{-1}DB$ acts
on $D\in \sD$ by cyclic permutation of diagonal blocks.
\end{proof}
\section{Proof of the results}\label{sec-thmain}
\subsection{Proof of Reduction~\protect\eqref{sreachtosreach}}
\label{sec-sreachtosreach}
Consider an instance of $\sreach(r,n)$, consisting
of a linear representation 
$(\alpha,\mu,\beta)$ of dimension $n$ over $\sS$,
together with $\gamma\in \sS$. 
Define the morphism $\emb\mu$ as in~\eqref{cons-blondel},
together with 
\begin{align}
\label{e-def-alphap}
\alpha'=(\alpha,\zero_{1n},\ldots,\zero_{1n})\in \sS^{1\times rn},
\;
\beta'=
\begin{pmatrix} \beta \\ \vdots \\
\beta 
\end{pmatrix}
\in \sS^{rn\times 1}
\enspace .
\end{align}
Then, it follows readily from Lemma~\ref{lemma-semidirect}
that 
\[
\alpha'\emb\mu(v)\beta'= 
\alpha \mu(a_1^{i_1} a_2^{i_2}\ldots a_{k+1}^{i_{k+1}})\beta
\enspace ,
\]
again with a cyclic indexing of $a_1,\ldots,a_r$.
Therefore, $\alpha'\emb\mu(v)\beta'$
takes the same values when $v\in \{b,c\}^*$
as $\alpha\mu(w)\beta$ when $w\in \Sigma_r^*$,
and $\sreachp(r,n)$ reduces 
to $\sreachp(2,rn)$. Using the equivalence~\eqref{equi-sreach},
we get that $\sreach(r,n)$ reduces to $\sreach(2,rn)$,
which shows~\eqref{sreachtosreach}.\qed
\subsection{Proof of Reduction~\eqref{vreachtovreach}}
\label{sec-vreachtovreach}
Consider an instance of $\vreach(r,n)$ consisting of a morphism
$\mu:\Sigma_r^*\to \sS^{n\times n}$ and vectors 
$\alpha,\eta\in\sS^{1\times n}$.
Define $\emb\mu$ as in~\eqref{cons-blondel},
$\alpha'$ as in~\eqref{e-def-alphap},
together with $\eta'=(\eta,\zero_{1n},\ldots,\zero_{1n})\in\sS^{1\times rn}$.
It follows from~\eqref{e-jj} that
\[ 
(\exists w\in \Sigma_r^*,\; \alpha\mu(w)= \eta)
\iff 
(\exists 0\leq k\leq r-1,\; \exists v\in \{b,c\}^*,\; 
\alpha'\emb\mu(v)=\eta' \emb\mu(b)^k )
\]
Thus, $\vreachp(r,n)$ reduces to $\vreachp(2,rn)$. 
Thanks to the equivalence~\eqref{equi-vreach}, this
shows that $\vreach(r,n)$ reduces to $\vreach(2,rn)$. 
\qed
\subsection{Proof of Reduction~\eqref{mreachtomreach}}
We consider an instance of $\mreach(r,n)$ consisting
of a morphism $\mu:\Sigma_r^*\to \sS^{n\times n}$
together with a matrix $M\in\sS^{n\times n}$. 
We shall split the proof in two cases.
\subsubsection{Case $M=\zero$}\label{subsec-btck}
Then, we apply the reduction
of~\cite[Th.~1]{blondel} and~\cite[Th.~1]{cassaigne},
which is valid over any semiring. For completeness,
we reprove this reduction. Consider the morphism
$\emb{\mu}:\Sigma_r^*\to \sS^{nr\times rn}$
built from $\mu$ as in \eqref{cons-blondel}. 
It follows readily from
Lemma~\ref{lemma-semidirect} that:
\begin{align}
(\exists w\in \Sigma_r^+,\;\mu(w)=\zero)
\iff
(\exists v\in \{b,c\}^+,\;\emb{\mu}(v)=\zero)
\enspace .
\end{align}
Therefore, $\mreach(r,n,M=\zero)$ reduces to $\mreach(2,rn,M=\zero)$.
\subsubsection{Case $M\neq\zero$}\label{subsec-ournew}
To any $1\leq i\leq r$, we associate the
morphism $\mui: \Sigma_r^*\to \sS^{n\times n}$
obtained from $\mu$ by exchanging the matrices
$\mu(a_i)$ and $\mu(a_r)$:
\[
\mui(a_r)=\mu(a_i),\; \mui(a_i)=\mu(a_r),\;
\mrm{ and }
\mui(a_j)=\mu(a_j),\mrm{ for }j\not\in\{i,r\} \enspace .
\]
We define the morphism $\eme\mui:\Sigma_r^*\to\sS^{rn\times rn}$
from $\mui$ as in~\eqref{eq-def-kappa}, and we set
\[
M'=\diag(M, \zero_{(r-1)n,(r-1)n})
\enspace .
\]
We claim that
\begin{align}
(\exists w\in \Sigma_r^+,\mu(w)=M)
\iff (\exists 1\leq i\leq r,\;\exists v\in \{b,c\}^+,\eme{\mui}(v)=M')
\enspace.
\label{ff}
\end{align}
Indeed, if $\mu(w)=M$ for some $w\in \Sigma_r^+$,
we write $w=a_iz$ with $z\in \Sigma_r^*$. Let $z'$
denote the word obtained from $z$ by exchanging
$a_i$ and $a_r$. Then, we get from~\eqref{e-news} that 
\[
\eme{\mui}\comp \delta(a_r z') 
= 
\begin{pmatrix} 
\mu(a_iz)& \zero_{n,(r-1)n}\\ 
\zero_{(r-1)n,n}& \zero_{(r-1)n,(r-1)n} 
\end{pmatrix} =M'\enspace,
\]
which shows the ``$\Rightarrow$'' implication
in~\eqref{ff}.
Conversely, let us assume that $\eme{\mui}(v)=M'$ 
for some $v\in\{b,c\}^+$. We can write (uniquely)
as in~\eqref{codefactor},
$v=c^{i_1}b\ldots c^{i_k}bc^{i_{k+1}}$, with $k\geq 0$.
Since $\eme{\mui}(c^r)=\zero$ and $M'\neq \zero$,
$v$ does not have $c^r$ as a factor, i.e.,
$i_1,\ldots,i_{k+1}\leq r-1$. 
If $i_{k+1}\neq 0$, 
using~\eqref{theformula} and the expression of $\eme{\mui}(c)$, 
we get
\[
\eme\mui(v)=
\eme\mui(v')\eme\mui(c^{i_{k+1}}) 
=\begin{pmatrix}
\zero_{nr,n} & * 
\end{pmatrix}
= M'
\]
and identifying the first diagonal block,
we get $\zero=M$, a contradiction.
Therefore, $k\geq 1$ and $v=c^{i_1}b\ldots c^{i_k}b=\delta(z)$,
where $z=a_{i_1+1}\ldots a_{i_k+1}$.
Then, using~\eqref{theformula} again,
we get that $\mui(z)=M$, hence $\mu(z')=M$,
where $z'$ is obtained from $z$ by exchanging
$a_i$ and $a_r$, which shows the ``$\Leftarrow$'' implication
in~\eqref{ff}. Then, the equivalence~\eqref{ff}
shows that $\mreach(r,n,M\neq\zero)$ reduces to $\mreach(2,rn,M\neq\zero)$.
\qed 
\subsection{Proof of Reduction~\eqref{vreachtomreach}}
Consider an instance of $\vreach(r,n)$,
which consists of $\alpha,\eta\in \sS^{1\times n}$
and a morphism $\mu:\Sigma_r^*\to \sS^{n\times n}$.
We associate to this instance the $(n+1)\times (n+1)$ matrix 
\begin{align}
M_\eta  = 
\begin{pmatrix} 
\zero_{11} & \eta \\  
\zero_{n1} & \zero_{nn} 
\end{pmatrix}, 
\end{align} 
and the morphism 
$\mu': \Sigma_{r+1}^*\to\sS^{(n+1)\times(n+1)}$,
defined by:
\begin{align}
\mu'(a_{r+1}) = 
\begin{pmatrix} 
\unit & \zero_{1n} \\  
\zero_{n1} & \zero_{nn} 
\end{pmatrix}, 
\;\mrm{and} \quad 
\mu'(a_i) = 
\begin{pmatrix} 
\zero_{11} & \alpha \mu(a_i) \\  
\zero_{n1} &\mu(a_i)   
\end{pmatrix} 
\quad \forall 1\leq i\leq r \enspace .
\end{align}
An immediate induction on the length of $w$ shows that 
\begin{align}
\label{e-meta}
\mu'(w)= 
\begin{pmatrix} 
\zero_{11} & \alpha \mu(w) \\  
\zero_{n1} &\mu(w)  
\end{pmatrix}, 
\;\mrm{and} \quad 
\mu'(a_{r+1}w) = 
\begin{pmatrix} 
\zero_{11}& \alpha \mu(w) \\  
\zero_{n1} & \zero_{nn}   
\end{pmatrix}
\quad \forall w\in\Sigma_{r}^+ \enspace .
\end{align} 
We claim that
if $\eta \neq \zero$, 
then
\begin{align}
\label{vmpfnzero}
(\exists w\in \Sigma_r^+,\;\alpha \mu(w)=\eta)
\iff (\exists z\in \Sigma_{r+1}^+,\;\mu'(z)=M_\eta) 
\enspace. 
\end{align}
The reduction $\vreach(r,n,\eta\neq\zero)\reducesto \mreach(r+1,n+1)$
will follow from~\eqref{vmpfnzero}.

Let us assume that $\alpha \mu(w)=\eta$ for some $w\in \Sigma_r^+$.
Then, it follows from~\eqref{e-meta} that $\mu'(z)=M_\eta$,
where $z=a_{r+1}w$, which shows
the ``$\Rightarrow$'' implication in~\eqref{vmpfnzero}.
Conversely, let us assume that $\mu'(z)=M_\eta $
for some $z\in \Sigma_{r+1}^+$. We can write 
$z=w_1a_{r+1}w_2a_{r+1}\ldots a_{r+1} w_{k+1}$,
where $w_1,\ldots, w_{k+1}\in \Sigma_r^*$ and $k\geq 0$.
Since
\[
\begin{pmatrix}
\zero & * \\
\zero & * 
\end{pmatrix}
\begin{pmatrix}
\unit &\zero \\
\zero &\zero 
\end{pmatrix}
= 
\zero\enspace,
\]
whatever the values of the ``$*$'' entries are,
and since $\mu'(z)=M_\eta\neq\zero$, it follows that if $k\geq 1$, 
then $w_1,\ldots , w_{k}$ must be equal to the empty word.
Therefore, since $\mu'(a_{r+1})^2=\mu'(a_{r+1})$, we can assume 
that $k\leq 1$. 
If $k=0$ we have 
$z=w_1$ and 
we readily check from~\eqref{e-meta}
that $\mu'(z)\neq M_\eta$, a contradiction.
Therefore $k=1$ and $z=a_{r+1}w$ for some $w\in \Sigma_r^*$. 
Since $\mu'(a_{r+1})\neq M_\eta$, it follows that
$z=a_{r+1}w$ for some $w\in \Sigma_r^+$. 
Then, we readily check from~\eqref{e-meta} that $\alpha\mu(w)=\eta $. 
This shows the 
``$\Leftarrow$'' implication
in~\eqref{vmpfnzero}. 

It remains to consider the case when $\eta=\zero$. 
Then, we introduce a trim unambiguous automaton $\sA$
recognizing $a_{r+1}\Sigma_r^*$.
We can take for $\sA$ the minimal automaton
of $a_{r+1}\Sigma_r^*$, which has two states,
$1$ and $2$, a set of initial
states $I=\{1\}$ and a set of final
states $F=\{2\}$,
with an associated morphism $\nu_{\sA}:\Sigma_{r+1}^*
\to\sS^{2\times 2}$ built as in~\S\ref{sec-cons-mat}:
\[
\nu_{\sA}(a_i)=
\begin{pmatrix}
\zero &\zero \\
\zero&\unit 
\end{pmatrix}
\mrm{ for }1\leq i\leq r
\mrm{ and }
\nu_{\sA}(a_{r+1})=
\begin{pmatrix}
\zero &\unit \\
\zero&\zero
\end{pmatrix},
\]
and 
\[
\sfM_{\sA}=\{M'\}
\mrm{ where }
M'= \begin{pmatrix}
\zero &\unit \\
\zero&\zero
\end{pmatrix} \enspace .
\]
Let $\mu''=\diag(\mu',\nu_{\sA})$ 
and $M''_\eta=\diag(M_\eta ,M')$.
By Lemma~\ref{lem-cons}, we have
\begin{align}
\mu''(z)=M''_\eta
\iff (z\in a_{r+1}\Sigma_r^* \mrm{ and }
\mu'(z)= M_\eta) \enspace .
\end{align}
But if $\mu'(z)=M_\eta = \zero_{nn}$ and $z=a_{r+1}w$
with $w\in \Sigma_r^*$, $w$ must be non-empty. Combining
this observation with~\eqref{e-meta}, we get
that
\[
(\exists z\in \Sigma_{r+1}^+,\;\mu''(z)=M''_\eta)
\iff (\exists w\in \Sigma_r^+,\;
\alpha \mu(w)= \eta =\zero ) \enspace ,
\]
which shows that $\vreach(r,n,\eta =\zero)\reducesto \mreach(r+1,n+3)$.
\qed
\subsection{Proof of Reduction~\eqref{sreachtovreach}}
Consider an instance of $\sreach(r,n)$ 
given by a series $s\in \se$ with a linear representation 
$(\alpha,\mu,\beta)$ of dimension $n$,
and $\gamma\in \sS$. 
By comparison to the proof of Reduction~\eqref{vreachtomreach},
we shall use a dual coding, and associate to this
instance the $1\times (n+1)$ matrices 
$\alpha'=(\alpha ,\zero_{11} )$ and $ \eta_\gamma =(\zero_{1n},\gamma)$, 
and the morphism 
$\mu': \Sigma_{r+1}^*\to\sS^{(n+1)\times(n+1)}$,
defined by:
\begin{align}
\mu'(a_{r+1}) = 
\begin{pmatrix} 
\zero_{nn} & \zero_{n1} \\  
\zero_{1n} & \unit 
\end{pmatrix}, 
\;\mrm{and} \quad 
\mu'(a_i) = 
\begin{pmatrix} 
\mu(a_i) & \mu(a_i)\beta \\  
\zero_{1n} & \zero_{11}  
\end{pmatrix} 
\quad \forall 1\leq i\leq r \enspace .
\end{align}
The dual version of~\eqref{e-meta} is:
\begin{align}
\label{e-vgamma}
\mu'(w)= 
\begin{pmatrix} 
\mu(w) & \mu(w)\beta \\  
\zero_{1n} & \zero_{11} 
\end{pmatrix}, 
\;\mrm{and} \quad 
\mu'(wa_{r+1}) = 
\begin{pmatrix} 
\zero_{nn} & \mu(w)\beta \\  
\zero_{1n} & \zero_{11}  
\end{pmatrix}
\quad \forall w\in\Sigma_{r}^+ \enspace .
\end{align} 
By dualizing the arguments of the proof of~\eqref{vmpfnzero},
we get that if $\gamma\neq\zero$, then
\begin{align}
\label{vpfnzero}
(\exists w\in \Sigma_r^+,\;s(w)=\gamma)
\iff (\exists z\in \Sigma_{r+1}^+,\;\alpha' \mu'(z)=\eta_\gamma) 
\enspace. 
\end{align}
The reduction $\sreach(r,n,\gamma\neq\zero)\reducesto \vreach(r+1,n+1)$
follows from~\eqref{vpfnzero}.

It remains to consider the case when $\gamma=\zero$. 
Then, we consider a trim unambiguous automaton $\sA$ 
with $2$ states recognizing $\Sigma_r^*a_{r+1}$, together
with the morphism $\nu_{\sA}:\Sigma_{r+1}^*\to \sS^{2\times 2}$
built as in \S\ref{sec-cons-mat}. 
We can assume that the initial state of $\sA$ is $1$ and that its final state is $2$. 
Let $\mu''=\diag(\mu',\nu_{\sA})$,
 $\alpha''=(\alpha,\unit ,\unit ,\zero$) 
and $\eta''=(\zero_{1n},\zero ,\zero ,\unit )$.
By Lemma~\ref{lem-cons}, 
\begin{align}
\alpha''\mu''(z)=\eta''
\iff (z\in \Sigma_r^* a_{r+1} \mrm{ and }
(\alpha ,\unit )\mu'(z)= (\zero_{1n},\zero ) ) \enspace .
\end{align}
But if $(\alpha ,\unit )\mu'(z)= (\zero_{1n},\zero )$ and $z=wa_{r+1}$
with $w\in \Sigma_r^*$, $w$ must be non-empty. Combining
this observation with~\eqref{e-vgamma}, we get
that
\[
(\exists z\in \Sigma_{r+1}^+,\;\alpha''\mu''(z)=\eta'')
\iff (\exists w\in \Sigma_r^+,\;
s(w)= \gamma =\zero ) \enspace ,
\]
which shows that $\sreach(r,n,\gamma=\zero)\reducesto \vreach(r+1,n+3)$.
\qed
\subsection{Proof of Reduction~\eqref{sreachtomreach}}
Consider an instance of $\sreach(r,n)$ 
given by a series $s\in \se$ with a linear representation 
$(\alpha,\mu,\beta)$ of dimension $n$,
and $\gamma\in \sS$. 
We associate to this instance the morphism 
$\mu': \Sigma_{r+1}^*\to\sS^{(n+2)\times(n+2)}$,
with $\mu'(a_i)$ as in~\eqref{eq-def-mup}, for $1\leq i\leq r$, 
and 
\[ 
\mu'(a_{r+1}) = 
\begin{pmatrix} \unit & \zero_{1n} & \zero_{11} \\  
\zero_{n1} & \zero_{nn} & \zero_{n1} \\  
\zero_{11} & \zero_{1n} & \unit 
\end{pmatrix} 
\enspace. \] 
Left and right multiplying~\eqref{eq-p-mup} by $\mu'(a_{r+1})$, 
we get: 
\begin{align}
\mu'(a_{r+1}wa_{r+1})= 
\begin{pmatrix} 
\zero_{11} & \zero_{1n}& s(w) \\  
\zero_{n1} & \zero_{nn} & \zero_{n1} \\  
\zero_{11} & \zero_{1n} & \zero_{11} 
\end{pmatrix}
\enspace,\qquad  \forall w\in\Sigma_{r}^+ \enspace. 
\label{e-mgamma}
\end{align}
Let
\[ M_{\gamma}=  
\begin{pmatrix} 
\zero_{11} & \zero_{1n} & \gamma \\  
\zero_{n1} & \zero_{nn} & \zero_{n1} \\  
\zero_{11} & \zero_{1n} & \zero_{11} 
\end{pmatrix} .\]
We claim that
if $\gamma\neq\zero$, 
then
\begin{align}
\label{pfnzero}
(\exists w\in \Sigma_r^+,\;s(w)=\gamma)
\iff (\exists z\in \Sigma_{r+1}^+,\; \mu'(z)=M_\gamma) 
\enspace. 
\end{align}
The reduction $\sreach(r,n,\gamma\neq\zero)\reducesto \mreach(r+1,n+2)$
will follow from~\eqref{pfnzero}.

Let us assume that $s(w)=\gamma$ for some $w\in \Sigma_r^+$.
Then, it follows from~\eqref{e-mgamma} that $\mu'(z)=M_\gamma$,
where $z=a_{r+1}wa_{r+1}$, which shows
the ``$\Rightarrow$'' implication in~\eqref{pfnzero}.
Conversely, let us assume that $\mu'(z)=M_\gamma$
for some $z\in \Sigma_{r+1}^+$. We can write 
\begin{align}
z=w_1a_{r+1}w_2a_{r+1}\ldots a_{r+1} w_{k+1}\enspace,
\label{fact}
\end{align}
where $w_1,\ldots, w_{k+1}\in \Sigma_r^*$ and $k\geq 0$.
Since $\mu'(a_{r+1})^2=\mu'(a_{r+1})$, 
 if some positive power $a_{r+1}^m$ appears
as a factor of $z$, we may replace this power by $a_{r+1}$
without changing $\mu'(z)$, which allows us to assume
that when $k\geq 2$, all the $w_2,\ldots,w_k$ are non-empty words. 
We remark that
\begin{align}\label{remarkthat}
\begin{pmatrix}
\zero &* & * \\
\zero &* & * \\
\zero &\zero &\zero
\end{pmatrix}
\begin{pmatrix}
\unit &\zero & \zero \\
\zero &\zero & \zero \\
\zero &\zero &\unit
\end{pmatrix}
\begin{pmatrix}
\zero &* & * \\
\zero &* & * \\
\zero &\zero &\zero
\end{pmatrix}
= 
\zero\enspace,
\end{align}
whatever the values of the ``$*$'' entries are.
It follows from \eqref{remarkthat} and from $\mu'(z)=M_\gamma\neq \zero$ that $z$ has no factor of the form
$wa_{r+1}w'$, with $w$, $w'\in \Sigma^+_r$. Therefore, in the factorization \eqref{fact},  $k\leq 2$ and 
at most one $w_i$ 
is different from the empty word. 
If $k\leq 1$, we have 
$z=w_1a_{r+1}$, or $z=a_{r+1}w_2$,
or $z=w_1$, and in all these cases, 
we readily check from~\eqref{eq-p-mup}
that $\mu'(z)\neq M_\gamma$, a contradiction.
Therefore $k=2$ and $z=a_{r+1}w_2a_{r+1}$ with $w_2\in \Sigma_r^+$.
Using~\eqref{e-mgamma}, we get $s(w_2)=\gamma$. 
This shows the ``$\Leftarrow$'' implication
in~\eqref{pfnzero}. 

It remains to consider the case when $\gamma=\zero$. 
Then, we consider a trim unambiguous automaton $\sA$ with $3$ states
recognizing $a_{r+1}\Sigma_r^*a_{r+1}$, together
with the morphism $\nu_{\sA}:\Sigma_{r+1}^*\to \sS^{3\times 3}$
and the set $\sfM_{\sA}\subset\sS^{3\times 3}$ built 
as in \S\ref{sec-cons-mat}. Let $\mu''=\diag(\mu',\nu_{\sA})$,
and $\sfM'=\set{\diag(M_\gamma,N)}{N\in \sfM_{\sA}}$. 
By Lemma~\ref{lem-cons}, 
\begin{align}
\mu''(z)\in \sfM'
\iff (z\in a_{r+1}\Sigma_r^* a_{r+1} \mrm{ and }
\mu'(z)= M_{\gamma}) \enspace .
\end{align}
But if $\mu'(z)=M_{\gamma}$ and $z=a_{r+1}wa_{r+1}$
with $w\in \Sigma_r^*$, $w$ must be non-empty. Combining
this observation with~\eqref{e-mgamma}, we get
that
\[
(\exists z\in \Sigma_{r+1}^+,\mu''(z)\in \sfM')
\iff (\exists w\in \Sigma_r^+,s(w)= \gamma =\zero ) \enspace ,
\]
which shows that $\sreach(r,n,\gamma=\zero)\reducesto \mreach(r+1,n+5)$.
\qed
\subsection{Proof of Reduction~\eqref{mreachtovreach}}\label{sec-blondel}
Consider an instance of $\mreach(r,n)$ consisting
of a morphism $\mu:\Sigma_r^*\to \sS^{n\times n}$
together with a matrix $M\in\sS^{n\times n}$. 
We associate to this instance the morphism $\mu': \Sigma_r^*\to\sS^{rn\times rn}$ defined by 
$\mu'=\diag (\mu ,\ldots ,\mu )$.   
Let $\vect $ be the matrix to vector operation that develops a square matrix into a row vector by taking its rows one by one. Then we have
\begin{align}\label{remarkblondel}
\vect (\mu (a_iw))=\vect (\mu (a_i))\mu'(w),\mrm{ for all }1\leq i\leq r\mrm{ and } w\in \Sigma_r^*\enspace .
\end{align}
It follows from~\eqref{remarkblondel} that
\[ 
(\exists w\in \Sigma_r^+,\; \mu(w)=M)
\iff 
(\exists 1\leq k\leq r,\; \exists v\in \Sigma_r^*,\; 
\vect(\mu (a_k))\mu'(v)=\vect(M))\enspace .
\]
Thus, $\mreach(r,n)$ reduces to $\vreachp(r,rn)$. 
Thanks to the equivalence~\eqref{equi-vreach}, this
shows that $\mreach(r,n)$ reduces to $\vreach(r,rn)$. 
\qed
\subsection{Proof of Corollary~\ref{cor-cassaigne}}
\label{sec-cor-cassaigne}
We have the following chain of reductions:
\[
\begin{array}{ccc}
\creach(\npost,3,\gamma=0)&\reducesto& \sreach(\npost,3,\gamma=0)\\
&&\downarrow\\
\creach(2,3\npost+2,\gamma=0) & \reducesfrom&\sreach(2,3\npost,\gamma=0)
\end{array}
\]
This follows from~\eqref{easy2},\eqref{sreachtosreach},
and~\eqref{easy} (and from the fact that the value of $\gamma$ is preserved
in~\eqref{easy2},\eqref{sreachtosreach}, and~\eqref{easy}).
Since $\creach(\npost,3,\gamma=0)$
is undecidable~\cite{manna}, it follows that $\creach(2,3\npost+2,\gamma=0) $,
the zero corner problem for $2$ matrices of dimension $3\npost+2$,
is undecidable.
\qed
\subsection{Proof of Proposition~\ref{cor-vmort}}\label{sec-proof-cor-vmort}
We shall combine a slight modification
of the proof of Paterson~\cite{paterson70}
with the idea of Bournez and Branicky~\cite{BB02} of using
the Modified Post Correspondence Problem.
Recall that the Modified Post Correspondence Problem (MPCP)
can be stated as:
given a finite set of pairs of words 
$\set{(u_i,v_i)}{1\leq i\leq r}$ over a finite alphabet,
is there a finite sequence $1\leq i_2,\ldots,i_k\leq r$
such that $u_{1}u_{i_2}\cdots u_{i_k}=v_{1}v_{i_2}\cdots v_{i_k}$? 
Of course, the MPCP is undecidable for any  value of $r$ for which
the Post Correspondence Problem is undecidable.
We shall assume, without loss of generality, that the
alphabet is $\Sigma=\{1,\ldots,n\}$. Let $b$ denote
any integer (strictly) greater than $n$, and
for any $w\in \Sigma^*$, let $[w]_b$ denote
the integer obtained by interpreting
the word $w$ in base $b$ (we set $[\_]_b=0$
in the case of the empty word), and let $|w|$ denote
the length of a word $w$. 
Paterson associated to any $u,v\in \Sigma^*$,
the matrix
\[
W(u,v)= 
\begin{pmatrix} 
b^{|u|} & 0 & 0 \\  
0 & b^{|v|} & 0 \\  {}
[u]_b & [v]_b & 1 
\end{pmatrix}
\]
and observed that
\begin{align}
\forall 
u',v'\in \Sigma^*
\enspace, 
\begin{pmatrix}
[u']_b & [v']_b & 1 
\end{pmatrix}
W(u,v)=
\begin{pmatrix}
[u'u]_b & [v'v]_b & 1 
\end{pmatrix} \enspace.
\label{e-paterson}
\end{align}
To any instance $I=\set{(u_i,v_i)}{1\leq i\leq r}$
of the MPCP over the alphabet $\Sigma$,
we associate the vector $\alpha=([u_1]_b,[v_1]_b,1)\in\Z^{1\times 3}$, 
and the morphism $\mu:\Sigma_{r+1}^*\to \Z^{3\times 3}$,
such that $\mu(a_i)=W(u_i,v_i)$ for $1\leq i\leq r$,
and $\mu(a_{r+1})=T$, where
\[
T= \begin{pmatrix}
1 & -1 &0\\
-1 & 1 & 0\\
0 & 0 & 0
\end{pmatrix}  \enspace.
\]
It follows readily from~\eqref{e-paterson}
and from the form of $T$ that
for all $1\leq i_2,\ldots,i_k\leq r$,
\begin{align}
\alpha W(u_{i_2},v_{i_2})
\cdots W(u_{i_k},v_{i_k})
T
=
\begin{pmatrix}
[u]_b-[v]_b & 
[v]_b-[u]_b & 
0 
\end{pmatrix}
\label{e-patnew}\\
\nonumber 
\mrm{where } 
u=u_1u_{i_2}\cdots u_{i_k}
\mrm{ and } v= v_1v_{i_2}\cdots v_{i_k}
\enspace . 
\end{align}
We claim that 
\begin{align}
\mrm{Instance }I\mrm{ has a solution}
\iff 
\exists w\in \Sigma_{r+1}^+,
\alpha\mu(w)=0
\enspace .
\label{reducpcp}
\end{align}
Indeed, the
``$\Rightarrow$'' implication in~\eqref{reducpcp}
follows readily from~\eqref{e-patnew}.
Conversely, let us assume that
$\alpha\mu(w)=0$ for some $w\in \Sigma_{r+1}^+$,
that we choose of minimal length. 
We can write $w=w_1a_{r+1}w_2\ldots
\ldots a_{r+1}w_{k+1}$, where $w_1,\ldots,w_{k+1}
\in \Sigma_{r}^*$. Since the matrices
$\mu(a_i)=W(u_i,v_i)$ all are invertible, 
for $1\leq i\leq r$, $\mu(w_{k+1})$ is invertible,
so that by minimality
of $w$, 
$w_{k+1}$ must be equal to the empty word.
We also note that all the matrices
$W(u_i,v_i)$ are of the form
\begin{align}
\begin{pmatrix} 
p & 0 & 0 \\  
0 & s & 0 \\  
q & t & 1 
\end{pmatrix}
\label{matrixUV}
\end{align}
where $p,s\geq 1$ and $q,t\geq 0$ 
and that the matrices of this form yield a semigroup.
Since 
\[
T
\begin{pmatrix} 
p & 0 & 0 \\  
0 & s & 0 \\  
q & t & 1 
\end{pmatrix}
T
=(p+s)T 
\]
we conclude, using again the minimality
of $w$, that $w_2$, \ldots, $w_k$ must be equal to
the empty word. Thus, $w=w_1a_{r+1}^m$, for some $m\geq 0$.
Since $T^2=2T$, the minimality
of $w$ yields $m\leq 1$.
If $m=0$, then, $w=w_1\in \Sigma_r^*$, and
$\alpha\mu(w)=(*,*,1)\neq 0$, a contradiction.
Thus $m=1$, so that $w=w_1a_{r+1}$
where $w_1=a_{i_2}\ldots a_{i_k}$ 
for some $1\leq i_2,\ldots,i_k\leq r$.
Then it follows from~\eqref{e-patnew} that 
$[u_1u_{i_2}\ldots u_{i_k}]_b=
[v_1v_{i_2}\ldots v_{i_k}]_b$,
which shows that $i_2,\ldots,i_k$ solves
Instance $I$ of the MPCP. We have proved
the ``$\Leftarrow$'' implication
in~\eqref{reducpcp}. 
\qed
\subsection{Proof of Corollaries~\protect\ref{cor-proj} and~\protect\ref{cor-proj2}}
\label{sec-proj}
Consider an instance of $\mreach(r,n)$ over $\zmax$,
consisting of a morphism $\mu: \Sigma_r^*\to \zmax^{n\times n}$
and a matrix $M\in \zmax^{n\times n}$. 
Let $\nu=\diag(\unit, \mu)$, and $M'=\diag(\unit,M)$.
Since
\[
\nu(w)\sim M'\iff \mu(w)=M \enspace,
\]
$\mreach(r,n)$ over $\zmax$,
which is undecidable when $r=2$ by Theorem~\ref{th-1},
reduces to the projective matrix reachability
problem for $r$ matrices of dimension $n+1$. This shows Corollary~\ref{cor-proj}. 

Consider now an instance of $\vreach(r,n)$
over $\zmax$, consisting of a morphism $\mu: \Sigma_r^*\to \zmax^{n\times n}$
and vectors $\alpha,\eta\in \zmax^{1\times n}$.
Let $\nu=\diag(\unit,\mu)$ as above,
$\alpha'=(\unit,\alpha)\in\zmax^{1\times (n+1)}$, and 
$\eta'=(\unit,\eta)\in\zmax^{1\times (n+1)}$.
Since
\[
\alpha'\nu(w)\sim \eta'\iff \alpha\mu(w)=\eta\enspace,
\]
$\vreach(r,n)$ over $\zmax$,
which is undecidable when $r=2$ by Theorem~\ref{th-1},
reduces to the projective vector reachability
problem for $r$ matrices of dimension $n+1$.
This shows Corollary~\ref{cor-proj2}. 
\subsection{Proof of Theorem~\protect\ref{th-4}} \label{sec-decidable}
We shall use the following stronger form of the
separation property.
\begin{lemma}\label{effective}
A semiring $\sS$ is (effectively) separated by morphisms
of finite image if, and only if, 
there are (effective) maps $B\mapsto \sS_B$
and $B\mapsto \pi_B$ which to any finite subset $B$ of $\sS$,
associate a finite semiring 
$\sS_B$ and a semiring morphism $\pi_B$ from $\sS$ to $\sS_B$,
such that $\pi_B^{-1}(\pi_B(y))=\{ y\}$ for all $y\in B$.
\end{lemma} 
\begin{proof}
The ``if'' part is trivial. 
Conversely, assume that a semiring $\sS$
is separated by morphisms of finite image,
and let $B=\{b_1,\ldots,b_k\}$ be a finite subset of $\sS$.
For all $1\leq i\leq k$, there is a finite semiring $\sS_{b_i}$ and a semiring morphism
$\pi_{b_i}$ from $\sS$ to $\sS_{b_i}$
such that $\pi_{b_i}^{-1}(\pi_{b_i}(b_i))=\{ b_i\}$. 
Let $\sS_B=\sS_{b_1}\times \cdots \times \sS_{b_k}$
denote the {\em Cartesian product} of $\sS_{b_1},\ldots,
\sS_{b_k}$, that is, the Cartesian product of the underlying sets
equipped with entrywise sum and product,
and consider the semiring morphism
$\pi_B: \sS\to \sS_B,\; \pi_B(y)=(\pi_{b_i}(y))_{1\leq i\leq k}$.
We have $\pi_B^{-1}(\pi_B(b_i))=\{ b_i\}$ for all $1\leq i\leq k$. 
This shows the ``only if'' part. Finally, we note that effective
aspects are preserved in the above construction.
\end{proof}
We now prove Theorem~\ref{th-4}. 
By Theorem~\ref{th-main},
it suffices to show that the matrix reachability problem
is decidable over $\sS$. 
Let us consider an instance of $\mreach(r,n)$ over $\sS$,
consisting of a morphism $\mu: \Sigma_r^*\to \sS^{n\times n}$
and a matrix $M\in \sS^{n\times n}$. 
Define $B=\{M_{ij} \mid 1\leq i,j\leq n\} $. We know that there is a finite semiring $\sS_B$ and a semiring morphism $\pi_B$ from $\sS$ to $\sS_B$
such that $\pi_B^{-1}(\pi_B(y))=\{ y\}$ for all $y\in B$. 
We extend $\pi_B $ to a map from $\sS^{n\times n}$ to $\sS_B^{n\times n}$ by making $\pi_B $ act on each entry.
We note that the problem: 
\begin{align}
\label{newp}
\exists w\in \Sigma_r^+, \pi_B\circ \mu(w)=\pi_B(M) \enspace ? 
\end{align}
is decidable. Indeed, using the effective part of Lemma~\ref{effective},
we can compute the matrices $\pi_B(X)\in \sS_{B}^{n\times n}$, for
$X\in \mu(\Sigma_r)\cup \{M\}$, and we know
the addition and multiplication tables of the (finite) semiring $\sS_B$.
Therefore, we can compute the finite semigroup
$\pi_B\circ \mu(\Sigma_{r}^+)$,
and we can test whether it contains $\pi_B(M)$.

We finally show that 
\begin{align}
\label{freduc}
(\exists w\in \Sigma_r^+, \mu(w)=M) \iff (\exists w\in \Sigma_r^+, \pi_B\circ \mu(w)=\pi_B(M)) 
\enspace .\end{align} 
Clearly, if $\exists w\in \Sigma_r^+$ such that $\mu(w)=M$, then
$\pi_B\circ \mu(w)=\pi_B(M)$.
 Conversely, assume that $\pi_B\circ \mu(w)=\pi_B(M)$ for some $w\in \Sigma_r^+$. 
Then, $(\mu(w))_{ij}\in \pi_B^{-1}(\pi_B(M_{ij}))=\{ M_{ij}\}$
 for all $1\leq i,j\leq n$ 
and therefore $\mu(w)=M$, which shows~\eqref{freduc}.
It follows readily from~\eqref{freduc}
that the matrix reachability problem is decidable over $\sS$.\qed 

Theorem~\ref{th-4} can be thought of
as an extension of Krob's~\cite[Proposition~2.2]{krobfatou}, which
shows that if $s\in\ser{\nmin}$ (resp. $s\in\ser{\nmax}$),
for all $\gamma\in \nmin$ (resp. $\gamma\in \nmax$),
$\set{w\in\Sigma_r^*}{s(w)=\gamma}$ is a constructible 
rational language. 
We can in fact restate Theorem~\ref{th-4} in the following more precise
way: 
\begin{theorem}\label{th-5}
Let $\sS$ denote a semiring  that is effectively separated by morphisms
of finite image, let $\alpha,\eta\in\sS^{1\times n}$,
$\mu:\Sigma_r^*\to\sS^{n\times n}$ a morphism,
$\beta\in\sS^{n\times 1}$,
$M\in\sS^{n\times n}$, and $\gamma\in\sS$.
Then, the following sets all are constructible
rational languages:
\begin{align}
&\set{w\in\Sigma_r^*}{\alpha\mu(w)\beta =\gamma}\enspace,\label{el1}\\
&\set{w\in\Sigma_r^*}{\alpha\mu(w) =\eta}\enspace,\label{el2}\\
&\set{w\in\Sigma_r^*}{\mu(w) =M}
\enspace .
\end{align}
\end{theorem}
\begin{proof}
It follows from the proof of~\eqref{freduc} that:
\begin{align}
\set{w\in\Sigma_r^*}{\mu(w) =M}=(\pi_B\comp \mu)^{-1} \{ \pi_B(M) \} 
\enspace .\label{freduc2}
\end{align}
Now, recall that the Kleene-Sch\"utzenberger theorem
shows that a language of $\Sigma_r^*$ is rational if and only if it 
can be written as $\kappa^{-1}(F)$, where
$\kappa$ is a morphism from $\Sigma_r^*$
to a finite monoid $P$, and $F$ is a subset of $P$.
Taking $P=\sS_B^{n\times n}$, $F=\{\pi_B(M)\}$,
and $\kappa=\pi_B\comp\mu$,
it follows from~\eqref{freduc2} that 
$\set{w\in\Sigma_r^*}{\mu(w) =M}$ is rational,
and this rational language is constructible since
$P,F$ and $\kappa$ can be effectively computed.
This argument can be readily adapted to 
the languages~\eqref{el1},\eqref{el2}.
For instance, in the case of~\eqref{el1},
we can take a finite semiring $\sS_\gamma$ together with
a morphism $\pi_\gamma:\sS\to\sS_\gamma$ such that
$\pi^{-1}_\gamma(\pi_\gamma(\gamma))=\{\gamma\}$,
and note that
$\set{w\in\Sigma_r^*}{\alpha\mu(w)\beta =\gamma}=\kappa^{-1}(F)$
where $P=\sS_\gamma^{n\times n}$,
$F=\set{U\in P}{\pi_\gamma(\alpha) U\pi_\gamma(\beta) =\pi_\gamma(\gamma)}$, and
$\kappa=\pi_\gamma\comp \mu$. The adaptation in the case
of~\eqref{el2} is similar.
\end{proof}
\subsection{Proof of Proposition~\protect\ref{prop-1}}
\label{sec-proof-prop-1}
Let $\sS$ be any of the semirings $\nmin $, $\nmax $, $\nmaxb $, $\sL$,
$\N$, $\nbar$, and let $\gamma \in \sS$ 
be arbitrary. 
If $\gamma $ is a natural number, let us denote by $n$ 
any natural number strictly greater than $\gamma $. 
If $\gamma $  is not a natural number choose $n$ arbitrarily 
(for example $n=1$). Consider the quotient of $\sS$ by the congruence 
which identifies 
all the integers greater than or equal to $n$.
(We call {\em congruence}
an equivalence relation which preserves the semiring structure.)
Let us denote by $\sS_\gamma$ the resulting finite semiring equipped
with the quotient laws, and by 
$\pi_\gamma $ the canonical morphism from $\sS$ to $\sS_\gamma$.
Then we have that $\pi_\gamma^{-1}(\pi_\gamma (\gamma ))=\{\gamma \}$,
which shows Proposition~\ref{prop-1}.
\qed
\subsection{Case $r=1$}\label{sec-r=1}
When $r=1$, the decidability of the reachability problems follows readily from known results. 
For instance, the cyclicity theorem for reducible max-plus matrices shows that if $A$ is a $n\times n$ matrix 
with entries in the semiring $\zmax$ there are positive integers $c,N$,
 such that for all 
$1\leq i,j\leq n$, there are scalars $\lambda_0,\ldots, \lambda_{c-1}$ 
(depending on $i$, $j$) such that 
for all $0\leq l\leq c-1$,
\begin{align} 
\forall n\geq N,\qquad (A^{(n+1)c+l})_{ij}=\lambda_l (A^{nc+l})_{ij} \enspace , \label{e-cycli} 
\end{align} 
and the integers $c,N$ together with the scalars $\lambda_l$ can be effectively computed. This cyclicity theorem, 
which is taken from~\cite[VI,1.1.10]{gaubert92a}, where it is proved more generally for matrices with entries in the semiring $(\Rm,\max,+)$, 
is an immediate consequence of the characterization of max-plus rational series in one letter as merge 
of ultimately geometric series, see~\cite{moller88}, \cite[VI,1.1.8]{gaubert92a} (or~\cite{gaubert94d}) and~\cite{krob93a}. 
It follows that the matrix, vector, and scalar reachability
problems in $\zmax$ are decidable when $r=1$. 
\newcommand{\etalchar}[1]{$^{#1}$}

\end{document}